\documentclass[12pt]{article}

\usepackage{amsmath, epsfig, cite}
\usepackage{amssymb}
\usepackage{amsfonts}
\usepackage{latexsym}

\newtheorem{thm}{Theorem}[section]

\newtheorem{conj}[thm]{Conjecture}
\newtheorem{lem}[thm]{Lemma}

\newcommand{\pf}{\noindent{\it Proof.} }

\numberwithin{equation}{section}

\newcommand{\qed}{{\hfill$\square$}\medskip}

\begin{document}

%\linenumbers

\begin{center}
{\Large\bf Some congruences involving powers of\\[5pt] Delannoy polynomials}
\end{center}

\vskip 2mm \centerline{Victor J. W. Guo}
\begin{center}
{\footnotesize Department of Mathematics, Shanghai Key Laboratory of
PMMP, East China Normal University,\\ 500 Dongchuan Road, Shanghai
200241,
 People's Republic of China\\
{\tt jwguo@math.ecnu.edu.cn,\quad
http://math.ecnu.edu.cn/\textasciitilde{jwguo}}}
\end{center}

%%date: November 18, 2006
%%date : June 29, 2002
%\vskip 5mm
%\noindent {\it Suggested Running title}: Two Identities of Gould
%\vskip 0.2cm \noindent{\it AMS Subject Classifications:} 05A10; 05A19

\vskip 0.7cm \noindent{\bf Abstract.} The Delannoy polynomial $D_n(x)$ is defined by
\begin{align*}
D_n(x)=\sum_{k=0}^{n}{n\choose k}{n+k\choose k}x^k.
\end{align*}
We prove that, if $x$ is an integer and $p$ is a prime not dividing $x(x+1)$, then
\begin{align*}
\sum_{k=0}^{p-1}(2k+1)D_k(x)^3 &\equiv  p\left(\frac{-4x-3}{p}\right) \pmod{p^2}, \\
\sum_{k=0}^{p-1}(2k+1)D_k(x)^4 &\equiv  p \pmod{p^2}, \\
\sum_{k=0}^{p-1}(-1)^k(2k+1)D_k(x)^3 &\equiv  p\left(\frac{4x+1}{p}\right) \pmod{p^2},
\end{align*}
where $\big(\frac{\cdot}{p}\big)$ denotes the Legendre symbol.
The first two congruences confirm a conjecture of Z.-W. Sun [Sci. China 57 (2014), 1375--1400]. The third congruence
confirms a special case of another conjecture of Z.-W. Sun [J. Number Theory 132 (2012), 2673--2699]. We also prove that,
for any integer $x$ and odd prime $p$, there holds
\begin{align*}
\sum_{k=0}^{p-1}(-1)^k(2k+1)D_k(x)^4 &\equiv
p\sum_{k=0}^{\frac{p-1}{2}} (-1)^k {2k\choose k}^2(x^2+x)^k(2x+1)^{2k} \pmod{p^2},
\end{align*}
and conjecture that it holds modulo $p^3$.


\vskip 2mm
\noindent{\it MR Subject Classifications}: 11A07, 11B65, 05A10

\section{Introduction}
The central Delannoy numbers (see \cite{CHV,Sulanke}) are defined by
\begin{align*}%\label{eq:delannoy}
D_n
=\sum_{k=0}^{n}{n+k\choose 2k}{2k\choose k}.
\end{align*}
Z.-W. Sun \cite{Sun0,Sun,Sun2}, among other things, proved many interesting congruences
on sums involving Delannoy numbers, such as
\begin{align*}
\sum_{k=0}^{p-1}(2k+1)D_k         &\equiv p+2p(2^{p-1}-1)-p(2^{p-1}-1)^2 \pmod{p^4}, \\
\sum_{k=0}^{n-1}(2k+1)D_k^2       &\equiv 0 \pmod{n^2},
\end{align*}
where $p$ is a prime greater than $3$. Z.-W. Sun \cite{Sun2} also introduced the Delannoy polynomial $D_n(x)$ as follows:
\begin{align*}
D_n(x)=\sum_{k=0}^{n}{n\choose k}{n+k\choose k}x^k,
\end{align*}
i.e., $D_n(x)=P_n(2x+1)$, where $P_n(x)$ is the Legendre polynomial of degree $n$ (see, for example, \cite[p.~1]{Koepf}).
Then he raised the following conjecture.
\begin{conj}{\rm\cite[Conjecture 5.1]{Sun2}} Let $x$ be an integer and let $m$ and $n$ be positive integers. Then
\begin{align}
\sum_{k=0}^{n-1}(2k+1)D_k(x)^m &\equiv  0 \pmod{n}. \label{eq:conj-sun-1}
\end{align}
If $p$ is a prime not dividing $x(x+1)$, then
\begin{align}
\sum_{k=0}^{p-1}(2k+1)D_k(x)^3 &\equiv  p\left(\frac{-4x-3}{p}\right) \pmod{p^2},  \label{eq:conj-sun-2} \\[5pt]
\sum_{k=0}^{p-1}(2k+1)D_k(x)^4 &\equiv  p \pmod{p^2},  \label{eq:conj-sun-3}
\end{align}
where $\big(\frac{\cdot}{p}\big)$ denotes the Legendre symbol.
\end{conj}

The congruence \eqref{eq:conj-sun-1} in a more general form has been confirmed by Pan \cite{Pan} recently.
However, Pan \cite{Pan} did not give an
integer coefficient polynomial formula for
$$
\frac{1}{n}\sum_{k=0}^{n-1}(2k+1)D_k(x)^m.
$$
In this paper we shall prove the following results.
\begin{thm}\label{thm:main1}
Let $n$ be a positive integer. Then
\begin{align}
&\hskip -1mm \frac{1}{n}\sum_{k=0}^{n-1}(2k+1)D_k(x)^3  \notag\\
&=\sum_{i=0}^{n-1} \sum_{j=0}^{n-1}\sum_{k=0}^{i}{n\choose j+k+1}{n+j+k\choose j+k}{i+j\choose i}{j\choose i-k}{j+k\choose k}
{2i\choose i}x^{i+j}(x+1)^{i},  \label{eq:sum-1} \\[5pt]
&\hskip -1mm \frac{1}{n}\sum_{k=0}^{n-1}(2k+1)D_k(x)^4  \notag\\
&=\sum_{i=0}^{n-1} \sum_{j=0}^{n-1}\sum_{k=0}^{i}{n\choose j+k+1}{n+j+k\choose j+k}{i+j\choose i}{j\choose i-k}{j+k\choose k}
{2i\choose i}{2j\choose j}(x^2+x)^{i+j}.  \label{eq:sum-2}
\end{align}
\end{thm}

\begin{thm}\label{thm:main2}
The supercongruences \eqref{eq:conj-sun-2} and \eqref{eq:conj-sun-3} are true.
\end{thm}

\begin{thm}\label{thm:main4}
Let $x$ be an integer and $p$ an odd prime. Then
\begin{align}
\sum_{k=0}^{p-1}(-1)^k(2k+1)D_k(x)^3 &\equiv  p\left(\frac{4x+1}{p}\right) \pmod{p^2},
\ \text{provided that $p\nmid x(x+1)$}, \label{eq:new-1} \\[5pt]
\sum_{k=0}^{p-1}(-1)^k(2k+1)D_k(x)^4 &\equiv
p\sum_{k=0}^{\frac{p-1}{2}}(-1)^k{2k\choose k}^2 (x^2+x)^k(2x+1)^{2k} \pmod{p^2}.  \label{eq:new-2}
\end{align}
\end{thm}

For any positive integer $n$ and $p$-adic integer $x$, Z.-W. Sun \cite[(4.6)]{Sun} conjectured that
\begin{align}
\nu_p\left(\frac{1}{n}\sum_{k=0}^{n-1}(-1)^k(2k+1)D_k(x)^3 \right)
\geqslant\min\{\nu_p(n),\nu_p(4x+1)\}, \label{eq:sun-last}
\end{align}
where $\nu_p(x)$ denotes the $p$-adic valuation of $x$. It is clear that the congruence \eqref{eq:new-1}
confirms  the $n=p$ case of \eqref{eq:sun-last}.

\section{Proof of Theorem \ref{thm:main1}}
It is easy to see that (see \cite[Lemma 3.2]{SunZH1.5})
\begin{align}
D_n(x)^2
=\sum_{k=0}^n{n\choose k}{n+k\choose k}{2k\choose k}x^k(x+1)^k, \label{eq:dnx-square}
\end{align}
which can be deduced from Clausen's formula \cite{Clausen} (with $a=-\frac{n}{2}$, $b=\frac{n+1}{2}$ and $x\to -4x(x+1)$):
\begin{equation}
_2F_1\left[\begin{array}{c} a,\,b \\ a+b+\frac{1}{2}\end{array};x\right]^2
={}_3F_2\left[\begin{array}{c} 2a,\,2b,\,a+b \\ 2a+2b,\,a+b+\frac{1}{2}\end{array};x\right]
,\quad |x|<1,  \label{eq:Clausen}
\end{equation}
and the following quadratic transformation of Gauss hypergeometric function (see \cite[p.~180]{Koepf}):
\begin{align}
_2F_1\left[\begin{array}{c} a,\,b \\ a+b+\frac{1}{2}\end{array};4x(1-x)\right]
={}_2F_1\left[\begin{array}{c} 2a,\,2b \\ a+b+\frac{1}{2}\end{array};x\right].  \label{eq:trans}
\end{align}

Writing  $D_\ell(x)^3$ as $D_\ell(x)^2\cdot D_\ell(x)$ and applying \eqref{eq:dnx-square}, we have
\begin{align}
&\hskip -2mm \frac{1}{n}\sum_{\ell=0}^{n-1}(2\ell+1)D_\ell(x)^3   \notag\\
&=\frac{1}{n}\sum_{\ell=0}^{n-1}(2\ell+1)\sum_{i=0}^\ell {\ell\choose i}{\ell+i\choose i}{2i\choose i}x^i(x+1)^i
\sum_{j=0}^\ell {\ell\choose j}{\ell+j\choose j}x^j.   \label{eq:proof-sum-1}
\end{align}
Note that (see the proof of \cite[Lemma 4.2]{GZ2})
\begin{align}
{\ell\choose i}{\ell+i\choose i}{\ell\choose j}{\ell+j\choose j}
=\sum_{k=0}^i {i+j\choose i}{j\choose i-k}{j+k\choose k}{\ell\choose j+k}{\ell+j+k\choose j+k}. \label{eq:from-GZ}
\end{align}
Moreover, by induction on $n$, we can easily prove that
\begin{align}
\sum_{\ell=k}^{n-1}(2\ell+1){\ell\choose k}{\ell+k\choose k} =n{n\choose k+1}{n+k\choose k}. \label{eq:proved-ind}
\end{align}
Substituting \eqref{eq:from-GZ} into \eqref{eq:proof-sum-1},
exchanging the summation order, and then utilizing \eqref{eq:proved-ind}, we complete the proof of \eqref{eq:sum-1}.

Similarly, writing  $D_\ell(x)^4$ as $D_\ell(x)^2\cdot D_\ell(x)^2$ and applying \eqref{eq:dnx-square},
we can prove \eqref{eq:sum-2}.

\section{Proof of Theorem \ref{thm:main2}}

\noindent{\it Proof of \eqref{eq:conj-sun-2}.} Letting $n=p$ be a prime in \eqref{eq:sum-1}, and noticing that
${p\choose k}\equiv 0\pmod p$ for $1\leqslant k\leqslant p-1$ and ${2p-1\choose p}\equiv 1\pmod p$, we obtain
\begin{align}
&\hskip -1mm \frac{1}{p}\sum_{k=0}^{p-1}(2k+1)D_k(x)^3  \notag\\
&=\sum_{i=0}^{p-1} \sum_{j=0}^{p-1}\sum_{k=0}^{i}{p\choose j+k+1}{p+j+k\choose j+k}{i+j\choose i}{j\choose i-k}{j+k\choose k}
{2i\choose i}x^{i+j}(x+1)^{i}  \notag \\
&\equiv \sum_{i=0}^{p-1} \sum_{j=0}^{p-1}{i+j\choose i}{j\choose p-i-1}{p-1\choose j}
{2i\choose i}x^{i+j}(x+1)^{i}  \pmod{p}.  \label{eq:1p-3}
\end{align}
For $0\leqslant i,j\leqslant p-1$, there holds
\begin{align}
{i+j\choose i}{j\choose p-i-1}
\begin{cases}=0, &\text{if $i+j<p-1$,} \\
\equiv 0\pmod{p}, &\text{if $i+j\geqslant p$.}
\end{cases} \label{eq:cases}
\end{align}
Therefore, the possible nonzero summands in \eqref{eq:1p-3} must satisfy $i+j=p-1$.
In other words, the congruence \eqref{eq:1p-3} may be simplified as
\begin{align*}
\frac{1}{p}\sum_{k=0}^{p-1}(2k+1)D_k(x)^3
&\equiv \sum_{i=0}^{p-1} {p-1\choose i}{p-1\choose p-i-1}{2i\choose i}x^{p-1}(x+1)^{i} \\
&\equiv \sum_{i=0}^{p-1} {2i\choose i}(x+1)^{i}  \pmod{p},
\end{align*}
where we used the fact ${p-1\choose i}\equiv (-1)^i \pmod{p}$ and Fermat's little theorem.
The proof then follows from the congruence
\begin{align}
\sum_{k=0}^{p-1}{2k\choose k}x^k\equiv
\left(\frac{1-4x}{p}\right) \pmod p  \label{eq:sun-tauraso}
\end{align}
due to Sun and Tauraso \cite[Theorem 1.1]{ST2}  (see also \cite[Lemma 2.1]{Sun2}).  \qed

\medskip
\noindent{\it Proof of \eqref{eq:conj-sun-3}.} Let $n=p$ be a prime in \eqref{eq:sum-2}. Similarly to the proof of \eqref{eq:conj-sun-2},
we have
\begin{align*}
&\hskip -3mm \frac{1}{p}\sum_{k=0}^{p-1}(2k+1)D_k(x)^4  \\
&\equiv \sum_{i=0}^{p-1} \sum_{j=0}^{p-1}{i+j\choose i}{j\choose p-i-1}{p-1\choose j}
{2i\choose i}{2j\choose j}(x^2+x)^{i+j}   \\
&\equiv \sum_{i=0}^{p-1} {2i\choose i}{2p-2i-2\choose p-i-1} \quad\text{(by \eqref{eq:cases} and Fermat's little theorem)} \\
&\equiv 1  \pmod{p},
\end{align*}
where in the last step we used the following fact
\begin{align}
{2i\choose i}\equiv 0\pmod p\quad\text{for}\quad \frac{p-1}{2}<i<p,  \label{eq:fact}
\end{align}
and ${p-1\choose \frac{p-1}{2}}^2\equiv 1\pmod p$.   \qed

\section{Proof of Theorem \ref{thm:main4}}

We need the following two lemmas.
\begin{lem}\label{thm:main3}
Let $n$ be a positive integer. Then
\begin{align}
&\hskip -1mm \frac{1}{n}\sum_{k=0}^{n-1}(-1)^{n-k-1}(2k+1)D_k(x)^3  \notag\\
&=\sum_{i=0}^{n-1} \sum_{j=0}^{n-1}\sum_{k=0}^{i}{n-1\choose j+k}{n+j+k\choose j+k}{i+j\choose i}{j\choose i-k}{j+k\choose k}
{2i\choose i}x^{i+j}(x+1)^{i},  \label{eq:sum-11} \\[5pt]
&\hskip -1mm \frac{1}{n}\sum_{k=0}^{n-1}(-1)^{n-k-1}(2k+1)D_k(x)^4  \notag\\
&=\sum_{i=0}^{n-1} \sum_{j=0}^{n-1}\sum_{k=0}^{i}{n-1\choose j+k}{n+j+k\choose j+k}{i+j\choose i}{j\choose i-k}{j+k\choose k}
{2i\choose i}{2j\choose j}(x^2+x)^{i+j}.  \label{eq:sum-12}
\end{align}
\end{lem}
\pf
It is exactly similar to the proof of Theorem \ref{thm:main1}. The difference is that we need to replace
\eqref{eq:proved-ind} by the following identity:
\begin{align*}
\sum_{\ell=k}^{n-1}(-1)^{n-\ell-1}(2\ell+1){\ell\choose k}{\ell+k\choose k} =n{n-1\choose k}{n+k\choose k},
\end{align*}
which can also be proved by induction on $n$.
\qed

\begin{lem}Let $n$ be a positive integer. Then
\begin{align}
\sum_{k=0}^n{n\choose k}{2k\choose k}{2n-2k\choose n-k}
=\sum_{k=0}^{n}{2k\choose k}^2 {k\choose n-k}(-4)^{n-k}.  \label{eq:zeil}
\end{align}
\end{lem}
\pf Applying Zeilberger's algorithm (see \cite{Koepf,PWZ}), we find that both sides of \eqref{eq:zeil}
satisfy the following recurrence relation:
\begin{align*}
(n+2)^2S(n+2)-4(3n^2+9n+7)S(n+1)+32(n+1)^2 S(n)=0.
\end{align*}
Noticing that they also have the same initial values $S(0)=1$ and $S(1)=4$, we complete the proof.  \qed

\medskip
\noindent{\it Proof of \eqref{eq:new-1}.} Letting $n=p$ be a prime not dividing $x(x+1)$ in \eqref{eq:sum-11}, we have
\begin{align}
&\hskip -1mm \frac{1}{p}\sum_{k=0}^{p-1}(-1)^{k}(2k+1)D_k(x)^3  \notag\\
&=\sum_{i=0}^{p-1} \sum_{j=0}^{p-1}\sum_{k=0}^{i}{p-1\choose j+k}{p+j+k\choose j+k}{i+j\choose i}{j\choose i-k}{j+k\choose k}
{2i\choose i}x^{i+j}(x+1)^{i}  \notag \\
&\equiv \sum_{i=0}^{p-1} \sum_{j=0}^{p-1}\sum_{k=0}^{i}(-1)^{j+k}{i+j\choose i}{j\choose i-k}{j+k\choose k}
{2i\choose i}x^{i+j}(x+1)^{i} \pmod{p},  \label{eq:new-new}
\end{align}
where we used the fact that, for $0\leqslant j,k\leqslant p-1$,
$$
{p-1\choose j+k}{p+j+k\choose j+k}{j+k\choose k}\equiv (-1)^{j+k}{j+k\choose k} \pmod{p}.
$$
By the Chu-Vandermonde summation formula, we get
\begin{align}
\sum_{k=0}^i (-1)^k {j\choose i-k}{j+k\choose k}=(-1)^i.  \label{eq:chu}
\end{align}
Substituting \eqref{eq:chu} into \eqref{eq:new-new} and using the binomial theorem, we obtain
\begin{align}
&\hskip -1mm \frac{1}{p}\sum_{k=0}^{p-1}(-1)^{k}(2k+1)D_k(x)^3  \notag \\
&\equiv \sum_{i=0}^{p-1} \sum_{j=0}^{p-1}(-1)^{i+j}{i+j\choose i}{2i\choose i}x^{i+j}(x+1)^{i} \notag\\
&= \sum_{i=0}^{p-1} \sum_{j=0}^{p-1}\sum_{k=0}^i (-1)^{i+j}{i+j\choose i}{2i\choose i}{i\choose k}x^{i+j+k} \notag \\
&\equiv \sum_{m=0}^{3p-3} x^m \sum_{i=0}^{\min\{p-1,m\}} (-1)^i {2i\choose i}
\sum_{j=0}^{\min\{p-1,m-i\}}(-1)^j {i+j\choose i}{i\choose m-i-j}
\pmod{p}.  \label{eq:new-new-2}
\end{align}

Note that, if $m-i\leqslant p-1$, then
\begin{align*}
\sum_{j=0}^{\min\{p-1,m-i\}}(-1)^j {i+j\choose i}{i\choose m-i-j}
=\sum_{j=0}^{m-i}(-1)^j {i+j\choose i}{i\choose m-i-j}=(-1)^{m-i};
\end{align*}
while if $m-i\geqslant p$, then for $0\leqslant i,j\leqslant p-1$, there holds
${i+j\choose i}{i\choose m-i-j}\equiv 0\pmod p$. Hence, we may simplify \eqref{eq:new-new-2} to
\begin{align}
\frac{1}{p}\sum_{k=0}^{p-1}(-1)^{k}(2k+1)D_k(x)^3
\equiv \sum_{m=0}^{p-1} (-x)^m \sum_{i=0}^{m} {2i\choose i}
+\sum_{m=p}^{2p-2} (-x)^m \sum_{i=m-p+1}^{p-1} {2i\choose i}\pmod{p}.
\label{eq:new-new-3}
\end{align}
By \eqref{eq:sun-tauraso} and Fermat's little theorem, we have
\begin{align}
\sum_{m=0}^{p-1} (-x)^m \sum_{i=0}^{m} {2i\choose i}
&=\sum_{i=0}^{p-1} {2i\choose i} \sum_{m=i}^{p-1}(-x)^m  \notag \\
&=\sum_{i=0}^{p-1} {2i\choose i}\frac{(-x)^i-(-x)^p}{1+x} \notag \\
&\equiv \frac{1}{1+x}\left(\frac{1+4x}{p}\right)+\frac{x}{1+x}\left(\frac{-3}{p}\right) \pmod{p}, \label{eq:thm4-1}
\end{align}
and
\begin{align}
\sum_{m=p}^{2p-2} (-x)^m \sum_{i=m-p+1}^{p-1} {2i\choose i}
&=(-x)^p\sum_{m=0}^{p-2} (-x)^m \sum_{i=m+1}^{p-1} {2i\choose i} \notag\\
&=(-x)^p\sum_{i=1}^{p-1} {2i\choose i} \sum_{m=0}^{i-1} (-x)^m  \notag\\
&=(-x)^{p} \sum_{i=0}^{p-1} {2i\choose i}\frac{1-(-x)^{i}}{1+x} \notag\\
&\equiv \frac{-x}{1+x}\left(\frac{-3}{p}\right)+\frac{x}{1+x}\left(\frac{1+4x}{p}\right) \pmod{p}, \label{eq:thm4-2}
\end{align}
Substituting \eqref{eq:thm4-1} and \eqref{eq:thm4-2} into \eqref{eq:new-new-3}, we complete the proof.
\qed

\medskip
\noindent{\it Proof of \eqref{eq:new-2}.} Let $n=p$ be a prime in \eqref{eq:sum-12}. Then similarly to \eqref{eq:new-new} we have

\begin{align}
\frac{1}{p}\sum_{k=0}^{p-1}(-1)^k(2k+1)D_k(x)^4
&\equiv \sum_{i=0}^{p-1}\sum_{j=0}^{p-1}(-1)^{i+j}{i+j\choose i}{2i\choose i}{2j\choose j}(x^2+x)^{i+j} \nonumber\\
&=\sum_{n=0}^{p-1}(-1)^n (x^2+x)^n \sum_{i=0}^n  {n\choose i}{2i\choose i}{2n-2i\choose n-i} \pmod{p}, \label{eq:dk4-4}
\end{align}
where we used the fact that ${i+j\choose i}\equiv 0 \pmod{p}$ for $0\leqslant i,j\leqslant p-1$ and $i+j\geqslant p$.

By \eqref{eq:zeil}, the right-hand side of \eqref{eq:dk4-4} is equal to
\begin{align*}
&\hskip -2mm
\sum_{n=0}^{p-1}(-1)^n (x^2+x)^n \sum_{k=0}^n  {2k\choose k}^2{k\choose n-k}(-4)^{n-k}  \\
&=\sum_{k=0}^{p-1}(-1)^k {2k\choose k}^2  \sum_{n=k}^{p-1}{k\choose n-k} 4^{n-k} (x^2+x)^n \\
&\equiv \sum_{k=0}^{\frac{p-1}{2}}(-1)^k{2k\choose k}^2  \sum_{n=k}^{2k}{k\choose n-k} 4^{n-k} (x^2+x)^n  \\
&=\sum_{k=0}^{\frac{p-1}{2}}(-1)^k{2k\choose k}^2  (x^2+x)^k (1+4x+4x^2)^{k} \pmod{p},
\end{align*}
where we used the congruence \eqref{eq:fact} and the binomial theorem.
This completes the proof.  \qed

\section{Two open problems}
Motivated by \eqref{eq:sun-last}, we raise the following conjecture:
\begin{conj}Let $n$ be a positive integer and $x$ a $p$-adic integer. Then
\begin{align}
\nu_p\left(\frac{1}{n}\sum_{k=0}^{n-1}(2k+1)D_k(x)^3 \right)
\geqslant\min\{\nu_p(n),\nu_p(4x+3)\}. \label{eq:guo-last}
\end{align}
\end{conj}
It is obvious that Theorem \ref{thm:main2} means that the $n=p$ case of \eqref{eq:guo-last} is true.

Finally, numerical calculation suggests the following refinement of \eqref{eq:new-2}.
\begin{conj}
Let $x$ be an integer and $p$ an odd prime. Then
\begin{align*}
\sum_{k=0}^{p-1}(-1)^k(2k+1)D_k(x)^4 \equiv
p\sum_{k=0}^{\frac{p-1}{2}}(-1)^k {2k\choose k}^2 (x^2+x)^k(2x+1)^{2k} \pmod{p^3}.
\end{align*}
\end{conj}

\vskip 5mm \noindent{\bf Acknowledgments.} This work was partially
supported by the Fundamental Research Funds for the Central
Universities and the National Natural Science Foundation of China
(grant 11371144).

\end{document}